\documentclass[12pt]{article}
\usepackage{amsthm, amsmath, amsfonts, amssymb, enumerate, enumitem}
\usepackage{tikz}
\usepackage{fullpage}
\usepackage[pagebackref]{hyperref}
\hypersetup{
    colorlinks=true,
    urlcolor=purple,
    linkcolor=purple,
    citecolor=purple
}

\usepackage{amsthm}

\theoremstyle{plain}
\newtheorem{thm}{Theorem}[section]          
\newtheorem{lem}[thm]{Lemma}

\newtheorem{con}[thm]{Conjecture}

\theoremstyle{definition}
\newtheorem{definition}[thm]{Definition}
\newtheorem{rem}[thm]{Remark}
\newtheorem{claim}{Claim}

\newcommand{\norm}[1]{\left\lVert #1\right\rVert}
\newcommand{\ip}[2]{\left\langle #1,#2\right\rangle}

\newcommand{\affl}[3]{%
\noindent #1,
\textsc{#3}\\
Email: \texttt{#2}\\[1.5pt]
}

\begin{document}
\title{$p$-energy lower bounds for vector chromatic numbers }
\author{Clive Elphick, Abhay Jayarajan and M. Rajesh Kannan  }

\maketitle
\begin{abstract}
Let $G$ be a simple graph on $n$ vertices and let $\mathcal{E}^+_p(G)$ and $\mathcal{E}^-_p(G)$ denote the positive and negative $p$-energy of $G$ respectively. Elphick, Tang and Zhang proved that for $p \ge 0$:

\[
\chi_f(G) \ge 1 +\max\left \{ \frac{\mathcal{E}_p^+(G)}{\mathcal{E}_p^-(G)}, \frac{\mathcal{E}_p^-(G)}{\mathcal{E}_p^+(G)} \right \}^{\frac{1}{|p - 1|}},
\]
where $\chi_f(G)$ denotes the fractional chromatic number and the right-hand side at $p=1$ is interpreted by continuity. They conjectured that this bound is valid for the vector chromatic number, $\chi_v(G)$, for $p \geq 2$. We prove the analogous bound for the strict vector chromatic number $\chi_{sv}(G)$ for every real $p \geq 2$. We also prove the conjecture for all integers $p \geq 2$ and reals $p \geq n-1$. We also extend the range of $p < 2$ for which counterexamples  are known. 
\end{abstract}



\section{Introduction}\label{sec:Introduction}

Let $G$ be a simple graph with $n$ vertices, $m$ edges and (unweighted) adjacency matrix $A$. Let $\lambda_1 \ge \cdots \ge \lambda_n$ denote the eigenvalues of $A$. The positive and negative $p$-energies of $G$  and the total $p$-energy of $G$ are defined as follows for $p \ge 0$:

\[
\mathcal{E}^+_p(G) = \sum_{\lambda_i  > 0}\lambda_i^p \mbox{  ;  } \mathcal{E}^-_p(G) = \sum_{\lambda_i < 0} |\lambda_i|^p \mbox{  and  } \mathcal{E}_p(G) = \mathcal{E}_p^+(G) + \mathcal{E}_p^-(G).
\]
When $p = 1$ the total energy corresponds to the much studied graph energy, (see \cite{Gutman2001} for a comprehensive survey) and when $p = 2$ the positive and negative energies correspond to positive and negative square energies, as studied for example in \cite{clive_miriam_felix_pawel_DAM_2016,AndoLin2015,pawel_clive_EJC_2023}. Recently higher order energies have received increasing attention. We refer to \cite{MR3141337, Nikiforov2016} for more details.

Several parameters lie between the clique number $\omega(G)$ and the chromatic number $\chi(G)$ of $G$; see Wocjan and Elphick~\cite{WocjanElphick2019} and Wocjan, Elphick and Anekstein~\cite{WocjanElphickAnek2023}. For the purposes of this paper we can simplify their analysis somewhat as follows:
\[
\omega(G) \le \chi_v(G) \le \chi_{sv}(G) \le \chi_{rv}(G) \le \xi_f(G) \le \chi_f(G) \le \chi(G),
\]
where $\chi_v$ is the vector chromatic number, $\chi_{sv}$ is the strict vector chromatic number, $\chi_{rv}$ is the rigid vector chromatic number, $\xi_f$ is the projective rank and $\chi_f$ is the fractional chromatic number. Let $\vartheta$ denote the Lov\'asz theta function (see \cite{Lovasz1979}). Note that $\chi_{sv}(G)$ equals the Lov\'asz theta function of the complement of $G$, that is $\vartheta(\overline{G})$ \cite{KargerMotwaniSudan1998}.

The projective rank of a graph is defined as follows.
\begin{definition}
    A $d/r$-representation of a graph $G$ is a collection of rank $r$, $d \times d$ projectors $\{P_v\}_{v \in V(G)}$ satisfying $P_u P_v=0$ for every edge $uv$. The projective rank of $G$ is defined by 
    $$\xi_f(G) = \inf \left \{ \dfrac{d}{r}: G \text{ has a } d/r-\text{representation }\right \}.$$
\end{definition}
Elphick, Tang and Zhang \cite[Theorem~1.8]{MR4971662} proved that for real $p \ge 0$:

\[
\xi_f(G) \ge 1 +\max\left \{ \frac{\mathcal{E}_p^+}{\mathcal{E}_p^-}, \frac{\mathcal{E}_p^-}{\mathcal{E}_p^+} \right \}^{\frac{1}{|p - 1|}}.
\]
At $p=1$, the right-hand side is defined by its limit as $p\to1$ \cite[Section~3.3]{MR4971662}. For the ordinary chromatic number, $p=0$ gives the inertia bound of Elphick and Wocjan \cite{ElphickWocjanInertia2017}, $p=2$ gives the square-energy bound of Ando and Lin \cite{AndoLin2015}, and $p\to\infty$ gives Hoffman's bound. 

The strict vector chromatic number of a graph is defined as follows.
\begin{definition}
    Let $k \in \mathbb{R}$ with $k \geq 2$. A strict vector $k$-coloring of a graph is a map that assigns to each vertex $u$ of $G$ a unit vector $t_u$ such that 
    $$\ip{t_u}{t_v} = \dfrac{-1}{k-1} \quad \quad \text{ for each edge }uv. $$
    The strict vector chromatic number of $G$ is defined by 
    $$\chi_{sv}(G)= \inf \left\{ k \geq 2 : G \text{ has a strict vector } k \text{-coloring } \right \}.$$
\end{definition}
The vector chromatic number of a graph is defined as follows.
\begin{definition}
    Let $k \in \mathbb{R}$ with $k \geq 2$. A vector $k$-coloring of a graph is a map that assigns to each vertex $u$ of $G$ a unit vector $t_u$ such that 
    $$\ip{t_u}{t_v} \leq \dfrac{-1}{k-1} \quad \quad \text{ for each edge }uv. $$
    The vector chromatic number of $G$ is defined by 
    $$\chi_v(G)= \inf \left\{ k \geq 2 : G \text{ has a vector } k \text{-coloring } \right \}.$$
\end{definition}

Elphick, Tang and Zhang also noted that $C_5$ demonstrates that the lower bound for the projective rank is  not a lower bound for $\chi_v(G)$ for $0 \leq p \leq 1.154$. However  they proposed the following conjecture \cite[Conjecture 6.3]{MR4971662}.

\begin{con}\label{con:chi_v_conj}
For all $p \ge 2$, the vector chromatic number of any graph $G$ satisfies
\[
\chi_v(G) \ge 1 +\max\left \{ \frac{\mathcal{E}_p^+(G)}{\mathcal{E}_p^-(G)}, \frac{\mathcal{E}_p^-(G)}{\mathcal{E}_p^+(G)} \right \}^{\frac{1}{p - 1}}.
\]
\end{con}

Coutinho, Spier and Zhang \cite[Theorem~1]{coutinho2024conicprogrammingunderstandsums}
proved the case \(p=2\) of this conjecture. As \(p\to\infty\),
the conjectured inequality reduces to Hoffman's bound for the
vector chromatic number,
\[
    \chi_v(G) \geq 1+\frac{\lambda_1}{|\lambda_n|},
\]
which was proved by Galtman \cite{Galtman2000} and,
independently, by Bilu \cite{Bilu2006}.

The organization of the paper is as follows. In Section~\ref{sec:sv}, we prove the $\chi_{sv}(G)$ analogue of Conjecture~\ref{con:chi_v_conj} for $p \geq 2$, and in Section~\ref{sec:vc} we prove Conjecture~\ref{con:chi_v_conj} for $\chi_v(G)$ for integer $p \ge 2$ and reals in $[n-1, \infty)$. In Section~\ref{sec:counter} we show that, for each $1<p<3/2$, there is a graph for which the extension to this exponent fails even for \(\chi_{sv}\). 

\section{Strict vector chromatic number}\label{sec:sv}

In this section, we prove the analogue of  Conjecture~\ref{con:chi_v_conj} for the strict vector chromatic number $\chi_{sv}(G)$. A vector $z \in  \mathbb{R}^d$ is nonnegative (resp. positive), denoted by $z \geq 0$ (resp. $z > 0$),  if all the entries of $z$ are nonnegative (resp. positive). 
For a nonnegative vector $z$ and a real number $r>0$, the notation $z^{[r]}$ denotes the componentwise power
$(z_i^r)_{i=1}^d$. A $d \times d$ nonnegative matrix $B$ is stochastic if all the row sums of $B$ are equal to $1$, and $B$ is said to be doubly stochastic if both $B$ and $B^T$ are stochastic.
\begin{lem}\label{lem:stoc_alpha}
    Let $P$ be a symmetric stochastic $d \times d$ matrix. Suppose that $x \in \mathbb{R}^d \setminus \{0\}$ and $0 < \alpha \leq 1$ satisfy $Px=-\alpha x$. Then for every real $p \geq 2$, 
    $$ \alpha ^{p-1} \norm{x^+}_p^p \leq \norm{x^-}_p^p  \mbox{  and  } \quad 
    \alpha ^{p-1} \norm{x^-}_p^p \leq \norm{x^+}_p^p$$
    where $x^{+}=(x_1^+,\dots,x_d^+)$ and $x^{-}=(x_1^-,\dots,x_d^-)$ with $x_i^+=\max\{x_i,0\}$ and $x_i^-=\max\{-x_i,0\}$.
\end{lem}
\begin{proof}
\setcounter{claim}{0}
    We prove the lemma in a series of claims.
    \begin{claim}\label{cl:Jensen}
        For $p>1$ and $y \in \mathbb{R}^d$, $\norm{Py}_p \leq \norm{y}_p$.
    \end{claim}
    \begin{proof}
        Since $p >1$, the function $u \mapsto |u|^p$ is a convex function. By Jensen's inequality we obtain
        \begin{align*}
            \norm{Py}_p^p&= \sum \limits_{i=1}^{d}\left | \sum \limits_{j=1}^{d} P_{ij}y_j \right|^p \\
            &\leq \sum \limits_{i=1}^{d} \sum \limits_{j=1}^{d} P_{ij} \left |y_j \right|^p \\
            &= \sum \limits_{j=1}^{d}  \left |y_j \right|^p \left( \sum \limits_{i=1}^{d} P_{ij} \right) \\
            &= \sum \limits_{j=1}^{d}  \left |y_j \right|^p= \norm{y}_p^p.
        \end{align*}
    \end{proof}
    Note that $x^+,x^- \geq 0 \text{ and } x=x^+-x^-.$
    Also note that,
    \begin{align*}
        -\alpha \sum \limits_{i=1}^d x_i&= \sum \limits_{i=1}^d \left(-\alpha x\right)_i \\
        &=\sum \limits_{i=1}^d \left(Px\right)_i = \sum \limits_{i=1}^d \sum \limits_{j=1}^d P_{ij}x_j\\
        &= \sum \limits_{j=1}^d \left( \sum \limits_{i=1}^d P_{ij} \right)x_j = \sum \limits_{j=1}^d x_j.
    \end{align*}
    Consequently, $\sum \limits_{i=1}^d x_i=0$. Hence $\norm{x^+}_p, \norm{x^-}_p >0.$

    \begin{claim}\label{cl:z_Pz_sum}
        For $p>1$ and $0 \leq z \in \mathbb{R}^d$, $\ip{z^{[p-1]}}{Pz}\leq \norm{z}_p^p.$
    \end{claim}
    \begin{proof}
        For $p>1$ put $q=p/(p-1)$. Since $z \geq 0$,
        $$\norm{z^{[p-1]}}_q= \left (\sum \limits_{i=1}^d z_i^{(p-1)q} \right) ^{\frac{1}{q}}=\left (\sum \limits_{i=1}^d z_i^{p} \right) ^{\frac{(p-1)}{p}}=\norm{z}_p^{p-1}.$$
        Therefore 
        $$\ip{z^{[p-1]}}{Pz} \leq \norm{z^{[p-1]}}_q \norm{Pz}_p \leq \norm{z}_p^{p-1}\norm{z}_p=\norm{z}_p^p.$$
        where the first inequality follows from H\"older's inequality and the second inequality follows from Claim~\ref{cl:Jensen}.
    \end{proof}
    Define $$\beta=\left( \frac{\norm{x^+}_p}{\norm{x^-}_p}\right)^{\frac{p}{p-1}}>0 \quad \quad \text{ and } \quad \overline{x}= \left(x^+ +\beta x^-\right)^{[p-1]}.$$
    By definition of $x^+, x^-$ it is easy to see that, $\overline{x}= \left(x^+\right)^{[p-1]}+\beta^{p-1} \left(x^-\right)^{[p-1]}$.
    \begin{claim}\label{cl:x-_lower}
        $$\ip{(x^-)^{[p-1]}}{Px^+} \leq \norm{x^-}_p^p.$$
    \end{claim}
    \begin{proof}
        We have
        \begin{align*}
            \ip{\overline{x}}{x}&= \ip{(x^+)^{[p-1]}+\beta^{p-1}(x^-)^{[p-1]}}{x^+-x^-}\\
            &=\ip{(x^+)^{[p-1]}}{x^+}-\beta^{p-1}\ip{(x^-)^{[p-1]}}{x^-} \\
            &= \norm{x^+}_p^p -\beta^{p-1} \norm{x^-}_p^p=0.
        \end{align*}
        Therefore, 
        $$\ip{\overline{x}}{Px^+}-\ip{\overline{x}}{Px^-}=\ip{\overline{x}}{P(x^+-x^-)}=\ip{\overline{x}}{Px}=-\alpha \ip{\overline{x}}{x}=0.$$
        That is, 
        \begin{equation}\label{eq:ip_P+_P-}
           \ip{\overline{x}}{Px^+}=\ip{\overline{x}}{Px^-}. 
        \end{equation}
        We also have,
        \begin{align*}
            \ip{\overline{x}}{Px^+}&= \ip{(x^+)^{[p-1]}+\beta^{p-1} (x^-)^{[p-1]}}{Px^+}\\
            &=\ip{(x^+)^{[p-1]}}{Px^+}+\beta^{p-1} \ip{(x^-)^{[p-1]}}{Px^+} \\
            &\geq \beta^{p-1} \ip{(x^-)^{[p-1]}}{Px^+}.
        \end{align*}
        Using Claim~\ref{cl:z_Pz_sum} and \eqref{eq:ip_P+_P-} we obtain,
        \begin{align*}
            \left (1+\beta \right) \beta^{p-1} \ip{(x^-)^{[p-1]}}{Px^+} &\leq \left(1+ \beta \right)\ip{\overline{x}}{Px^+}\\
            &= \ip{\overline{x}}{Px^+}+ \beta\ip{\overline{x}}{Px^-}\\
            &= \ip{ \left(x^+ +\beta x^-\right)^{[p-1]}}{P(x^+ + \beta x^-)} \\
            &\leq \norm{x^+ + \beta x^-}_p^p\\
            &= \norm{x^+}_p^p + \beta^p \norm{x^-}_p^p \\
            &=\beta^{p-1} \norm{x^-}_p^p + \beta^p \norm{x^-}_p^p \\
            &= (1+\beta)\beta^{p-1} \norm{x^-}_p^{p}.
        \end{align*}
        The claim follows by dividing by $(1+\beta)\beta^{p-1} > 0.$
    \end{proof}
    \begin{claim}\label{cl:x+_upper}
        Let $p \geq 2$. Then $\alpha^{p-1}\norm{x^+}_p^p \leq \ip{(x^-)^{[p-1]}}{Px^+}.$
    \end{claim}
    \begin{proof}
        If $x^+_i>0$, then $x^-_i=0$. Thus $Px^+-Px^-=-\alpha x^+ +\alpha x^-$ gives
        $$(Px^-)_i=(Px^+)_i+\alpha x^+_i \geq \alpha x^+_i.$$
        Also if $x^+_i=0$, then $(Px^-)_i \geq 0=\alpha x^+_i$. \\
        Therefore for each $1 \leq i \leq d$, we have, 
        \begin{equation}\label{eq:alphax_ileqP-}
            \alpha x^+_i \leq (Px^-)_i.
        \end{equation}
        Further on $[0, \infty)$, and with $p \geq 2$, the function $u \mapsto u^{p-1}$ is convex. Therefore by Jensen's inequality we obtain
        \begin{align}\label{eq:P-_Jensen}
            \left( (Px^-)_i \right)^{p-1} &= \left ( \sum \limits_{j=1}^{d} P_{ij}x_j^- \right ) ^{p-1} \notag \\
            &\leq \sum \limits_{j=1}^{d} P_{ij}(x_{j}^{-})^{p-1}.  \notag \\
        \end{align}
        Since $p-1>0$, the function $u \mapsto u^{p-1}$ is increasing on $[0, \infty)$. Therefore using \eqref{eq:alphax_ileqP-} and $\eqref{eq:P-_Jensen}$ we get
        $$(\alpha x_i^+)^{p-1} \leq \left ( (Px^-)_i \right) ^{p-1} \leq \sum \limits_{j=1}^{d} P_{ij}(x_j^-)^{p-1}.$$
        Multiplying the above inequality by $x^+_i$ and then summing over all indices gives
        \begin{align*}
            \alpha ^{p-1}\sum \limits_{i=1}^{d}(x^+_i)^p &\leq \sum \limits_{i=1}^{d}x^+_i \sum \limits_{j=1}^{d} P_{ij}(x_j^-)^{p-1} \\
            &= \sum \limits_{j=1}^{d} (x_j^-)^{p-1} \left ( \sum \limits_{i=1}^{d}P_{ji}x^+_i \right)\\
            &= \sum \limits_{j=1}^{d} (x_j^-)^{p-1}(Px^+)_j \\
            &= \ip{(x^-)^{[p-1]}}{Px^+}.
        \end{align*}
    \end{proof}
    Thus combining Claim~\ref{cl:x-_lower} and Claim~\ref{cl:x+_upper}, we get 
    $$\alpha^{p-1}\norm{x^+}_p^p \leq \norm{x^-}_p^p.$$

Since $P(-x)=-\alpha(-x)$ and $(-x)^+=x^-, (-x)^-=x^+,$ applying the established inequality to $-x$ gives $\alpha^{p-1}\|x^-\|_p^p\le\|x^+\|_p^p.$
\end{proof}

For a real square matrix $B$, by $B \succeq 0$ we mean $B$ is a symmetric positive semidefinite matrix. For two matrices $B$ and $C$ of the same order, $B \circ C$ denotes the entrywise product of the matrices $B$ and $C$. That is, $(B\circ C)_{ij} = b_{ij}\cdot c_{ij}. $\\
We now state and prove the main result in this section. 
\begin{thm}\label{thm:chi_sv}
    Let $G$ be a graph on $n \geq 3$ vertices and $m \geq 1$ edges with (unweighted) adjacency matrix $A$ and strict vector chromatic number $\chi_{sv}(G)$. Let $p \geq 2$. Then
    $$ \chi_{sv}(G) = \vartheta(\overline{G}) \geq 1 +\max\left \{ \frac{\mathcal{E}_p^+}{\mathcal{E}_p^-}, \frac{\mathcal{E}_p^-}{\mathcal{E}_p^+} \right \}^{\frac{1}{p - 1}}. $$
\end{thm}
\begin{proof}
    Let $\chi_{sv}(G) = k$.  By the Gram matrix formulation of strict vector chromatic number, there exists a Gram matrix $Q$ with $Q \succeq 0, Q_{ii}=1  \quad (1 \leq i \leq n), Q_{ij} = \frac{-1}{k-1}$ on every edge $ij$. Consequently, 
     \begin{equation*}
         Q \circ A= \frac{-1}{k-1}A.
     \end{equation*}
     Let $z^1, \dots, z^n$ be an orthonormal basis of $ \mathbb{R}^n$ such that $Az^i= \lambda_i z^i$. Define $P=(P_{ij})$ by
     \begin{equation*}
         P_{ij}=(z^i \circ z^j)^\top Q (z^i \circ z^j).
     \end{equation*}
     It is easy to see that
$$P_{ij} =\sum \limits_{l,t}Q_{lt}z^i_lz^j_l z^i_t z^j_t.$$
Thus clearly $P$ is symmetric and since $Q \succeq 0$, we have $P_{ij} \geq 0.$ Next consider,
\begin{align*}
    \sum \limits_{j} P_{ij} &= \sum \limits_{j}\sum \limits_{l,t}Q_{lt}z^i_lz^j_l z^i_t z^j_t \\
    &= \sum \limits_{l,t}Q_{lt}z^i_lz^i_t \sum \limits_{j}z^j_lz^j_t \\
    &= \sum \limits_{l,t}Q_{lt}z^i_lz^i_t \delta_{lt}  \\
    &= \sum \limits_{l} Q_{ll} z^i_lz^i_l =1.
\end{align*}
Thus, $P$ is a symmetric entrywise nonnegative doubly stochastic matrix. Next define the vector $\overline \lambda=(\lambda_1, \dots, \lambda_n) \in \mathbb{R}^n$. Then, by the spectral decomposition, we have $A=\sum \limits_{j=1}^{n} \lambda_j z^j (z^j)^{\top}$. Thus we have for each $i$,
\begin{align*}
    (P \overline\lambda)_i &= \sum \limits_{j=1}^n P_{ij} \lambda_j \\
    &= \sum \limits_{j=1}^{n} \left(\sum \limits_{l,t} Q_{lt} z^i_l z^j_l z^i_t z^j_t \right)\lambda_j \\
    &=\sum \limits_{l,t} Q_{lt} z^i_l z^i_t \sum \limits_{j=1}^{n} \lambda_j z^j_l z^j_t \\
     &=\sum \limits_{l,t} Q_{lt} z^i_l z^i_t A_{lt} \\
     &= (z^i)^\top (Q \circ A) z^i \\
     &=\frac{-1}{k-1} (z^i)^\top A z^i \\
     &= \frac{-1}{k-1} \lambda_i.
\end{align*}
Thus $P \overline \lambda = -\left(\frac{1}{\chi_{sv}(G)-1}\right) \overline \lambda$.
    By Lemma~\ref{lem:stoc_alpha}, for $p \geq 2$,
    $\left(\frac{1}{\chi_{sv}(G)-1}\right)^{p-1} \sum \limits_{i=1}^{n} (\lambda_i^+)^p \leq \sum \limits_{i=1}^{n} (\lambda_i^-)^p$. That is,
    $$ \left(\frac{1}{\chi_{sv}(G)-1}\right)^{p-1} \mathcal{E}_p^+(G) \leq \mathcal{E}_p^-(G).$$
    Therefore, 
    $$ \frac{\mathcal{E}_p^+(G)}{\mathcal{E}_p^-(G)} \leq (\chi_{sv}(G)-1)^{p-1}. $$
    Similarly, we get, \\
    $$ \frac{\mathcal{E}_p^-(G)}{\mathcal{E}_p^+(G)} \leq (\chi_{sv}(G)-1)^{p-1}. $$
    Consequently, 
    $$ 1 +\max\left \{ \frac{\mathcal{E}_p^+}{\mathcal{E}_p^-}, \frac{\mathcal{E}_p^-}{\mathcal{E}_p^+} \right \}^{\frac{1}{p - 1}} \leq \chi_{sv}(G).$$
\end{proof}
\begin{rem}

In the proof of Theorem~\ref{thm:chi_sv}, we use the doubly stochastic matrix constructed in the proof of \cite[Lemma~4]{WocjanElphickAnek2023}, indexed with respect to a common orthonormal eigenbasis. In the present setting, this matrix is also symmetric.
\end{rem}

\begin{rem}
Theorem~\ref{thm:chi_sv} also holds with $A$ replaced by any nonzero Hermitian weighted adjacency matrix $B$ of $G$. With $Q$ and $k$ as in the proof, choose an orthonormal eigenbasis $z^1,\ldots,z^n$ of $B$ in $\mathbb C^n$
and define
\[
P_{ij}=(z^i\circ\overline{z^j})^*Q(z^i\circ\overline{z^j}).
\]
Then $P$ is again symmetric and doubly stochastic, and $P \overline \lambda=-\overline \lambda/(k-1)$, where $\overline \lambda$ is the corresponding vector of eigenvalues of $B$. Lemma~\ref{lem:stoc_alpha} therefore gives the same bound with  the $p$-energies of $B$ for every $p\geq 2$.
\end{rem}
\section{Vector chromatic number}\label{sec:vc}

Let \(K\) be a real symmetric matrix of order \(d\), and let
\(w_1,\ldots,w_d\) be an orthonormal eigenbasis with corresponding
eigenvalues \(\mu_1,\ldots,\mu_d\). For a real-valued function \(f\)
on the spectrum of \(K\), define
\[
f(K)=\sum_{i=1}^{d}f(\mu_i)w_iw_i^\top.
\]
Write \(K=K_+-K_-\), where $K_+$ and $K_-$ are the positive and negative part of the spectral decomposition of $K$. Then $K_+,K_-\succeq0$ and \(K_+K_-=0\).
The real powers of \(K_+\) and \(K_-\) are therefore defined for
every \(p>0\). In particular, for the adjacency matrix \(A\) of \(G\),
\[
\mathcal E_p^+(G)=\operatorname{tr}\left((A_+)^p \right),
\qquad
\mathcal E_p^-(G)=\operatorname{tr}\left( (A_-)^p\right ).
\]

We are now ready to state the main result of the section. In view of the above identities we state the result in terms of traces of the positive and negative parts of the adjacency matrix of graph $G$. 
\begin{thm}\label{thm:partial}
    If $G$ is a graph on $n \geq 3$ vertices and $m \geq 1$ edges with unweighted adjacency matrix $A$ and vector chromatic number $\chi_v(G)$ then 
    $$\chi_v(G) \geq 1+ \left(\max \left \{ \dfrac{\operatorname{tr}\left((A_+)^p\right)}{\operatorname{tr}\left((A_-)^p\right)}, \dfrac{\operatorname{tr}\left((A_-)^p\right)}{\operatorname{tr}\left((A_+)^p\right)}\right \}\right)^{\frac{1}{p-1}}$$
    for $p\in\{2,3,\ldots\}\cup[n-1,\infty)$. 
\end{thm}

 Let \(\mathcal W\) consist of
all real symmetric matrices \(T=(t_{ij})\) such that
\[
 t_{ij}\ge0,\qquad
 t_{ij}=0\ \text{if }i=j\text{ or }\{i,j\}\notin E(G),\qquad
 \sum_{\{u,v\}\in E(G)}t_{uv}=1.
\]

Fix $p\in\{2,3,\ldots\}\cup[n-1,\infty)$. Note that, $\mathcal{W}$ is a compact set. For a real symmetric matrix $T$ having both positive and negative eigenvalues, write
\[
\kappa_T(p)=\max\left\{\frac{\operatorname{tr}\left((T_+)^p\right)}{\operatorname{tr}\left((T_-)^p\right)},\frac{\operatorname{tr}\left((T_-)^p\right)}{\operatorname{tr}\left((T_+)^p\right)}\right\}.
\]
Each $T\in\mathcal W$ is nonzero and traceless, so both traces in $\kappa_T(p)$ are positive and $T\mapsto \kappa_T(p)$ is continuous on the compact set $\mathcal W$. Let $W$ be a maximizer, with positive and negative parts $W_+$ and $W_-$. Then either $\kappa_1=\operatorname{tr}\left((W_+)^p\right)/\operatorname{tr}\left((W_-)^p\right)\ge1$ or $\kappa_2=\operatorname{tr}\left((W_-)^p\right)/\operatorname{tr}\left((W_+)^p\right)\ge1$. If $\kappa_1=\kappa_2=1$, then $\kappa_A(p)=\kappa_{A/m}(p)\le1$, and the theorem reduces to $\chi_v(G)\ge2$, which holds since $G$ has an edge. To treat the cases $\kappa_1>1$ and $\kappa_2>1$ together, put $D=W$ if $\kappa_1>1$ and $D=-W$ if $\kappa_2>1$. Thus $D=\sigma W$ for some $\sigma\in\{1,-1\}$, and with $\gamma=\operatorname{tr}\left((D_+)^p\right)$ and $\delta=\operatorname{tr}\left((D_-)^p\right)$ the ratio $\kappa=\gamma/\delta>1$ equals $\max \limits_{T\in\mathcal W}\kappa_T(p)$. Set $$F=(D_+)^{p-1}+\kappa\,(D_-)^{p-1}.$$

Throughout Lemmas~\ref{lem:F_0for_posW}--\ref{lem:L_upperlower}, the exponent $p\in\{2,3,\dots\}\cup[n-1,\infty)$ is the one fixed above, and $W$, $D=\sigma W$, $\gamma$, $\delta$, $\kappa$ and $F$ are as just defined.

\begin{lem}\label{lem:F_0for_posW}
If \(W_{ij}>0\), then \(F_{ij}=0\).
\end{lem}
\begin{proof}
For symmetric $X$, let $\varphi(X)=\operatorname{tr}\left((X_+)^p\right)-\kappa\operatorname{tr}\left((X_-)^p\right)$. Note that the functions $t\mapsto(t^+)^p$ and $t\mapsto(t^-)^p$ are continuously differentiable for $p>1$, with derivatives $p(t^+)^{p-1}$ and $-p(t^-)^{p-1}$. Hence the function $\varphi$ is differentiable, and for every real symmetric matrix $T$
\[
\frac{d}{d\varepsilon}\varphi(X+\varepsilon T)\Big|_{\varepsilon=0}=p\operatorname{tr}\big(\left((X_+)^{p-1}+\kappa (X_-)^{p-1}\right)T\big).
\]
 By the maximality of $\kappa$, $\varphi(\sigma V)\le0$ for every $V \in \mathcal{W}$, and $\varphi(D)=0$. Let $W_{ij}>0$. Define $T=e_ie_j^{T}+e_je_i^{T}$. 

For $|\varepsilon|<W_{ij}$, put
\[
V_\varepsilon=\frac{W+\varepsilon\sigma T}{1+\varepsilon\sigma}
\in\mathcal W.
\]
By positive homogeneity of $\varphi$ and maximality of $\kappa$,
\[
g(\varepsilon):=\varphi(D+\varepsilon T)
=(1+\varepsilon\sigma)^p\varphi(\sigma V_\varepsilon)
\le0=g(0).
\]
Hence $0=g'(0)=2p\cdot F_{ij}$.
\end{proof}

A real symmetric matrix is doubly nonnegative (DNN) if it is positive semidefinite and entrywise nonnegative. We recall the notion of divided differences used in the following results. 
\begin{definition}\cite[Section~1]{KomisarskiRajba2021}
Let $I$ be an interval and $f:I\to\mathbb R$. For distinct points $x_0,\ldots,x_k\in I$, define
\begin{equation*}
[x_0;f]=f(x_0)
\end{equation*}
and, for every integer \(k\ge1\),
\begin{equation*}
[x_0,\ldots,x_k;f]
=
\frac{[x_1,\ldots,x_k;f]-[x_0,\ldots,x_{k-1};f]}
     {x_k-x_0}.
\end{equation*}
The quantity $[x_0,\ldots,x_k;f]$ is a divided difference of
order $k$.
\end{definition}
\begin{lem}\label{lem:conv_divdiff}\cite[Proposition~2]{KomisarskiRajba2021}
    Let $I$ be an open interval and $k \geq 2$ be an integer. All divided differences of $f: I \rightarrow \mathbb{R}$ in $I$ of order $k$ are nonnegative if and only if $f^{(k-2)}$ exists and is convex on $I$.
\end{lem}
\begin{thm}{\cite[Theorem~2.1]{GuillotKhareRaja2013}}\label{thm:DNN_divdiff}
    Let $d \geq 2$ and $f(t) \geq 0$ be  continuous for $t \geq 0$. Then  $f(P)$ is DNN  for every DNN matrix $P$ of order $d$ if and only if all divided differences of $f$ on $[0, \infty)$ of order $1, 2, \cdots, d-1$ are nonnegative. 
\end{thm}

\begin{lem}\label{lem:DNNpower}
   Let $P$ be a real symmetric entrywise nonnegative matrix of order $d \geq 3$. Then $(P_+)^\alpha$ is entrywise nonnegative for every real $\alpha \geq d-2$.
\end{lem}
\begin{proof}
    Let $\lambda_{\min}(P)$ be the least eigenvalue of $P$. Choose $s \geq \max\{0, -\lambda_{\min}(P) \}$. Then $B= P+sI$ is DNN. Define the function $f$ by $t \mapsto \left( (t-s)_+\right)^{\alpha}$. \\
    The first divided differences of $f$ are nonnegative since $f$ is nondecreasing. For $2\leq k\leq d-1$,
\begin{equation*}
f^{(k-2)}(t)
=
\left(\prod_{j=0}^{k-3}(\alpha-j)\right)
(t-s)_+^{\alpha-k+2},
\end{equation*}
where an empty product equals $1$. The coefficient is positive, and $\alpha-k+2\geq 1$. Hence $f^{(k-2)}$ exists and is convex on $\mathbb R$.
    By Lemma~\ref{lem:conv_divdiff} the divided differences of $f$ of order $k$ are nonnegative. Thus, by Theorem~\ref{thm:DNN_divdiff} $f(B)$ is DNN. In particular, $f(B)= (P_+)^{\alpha}$ is entrywise nonnegative.
\end{proof}

\begin{lem}\label{lem:D_+D_geq0}
    For each $i$ and $j$, the following inequality is true
    $$(\left(D_+ \right)^{p-1})_{ij}D_{ij} \geq 0.$$
\end{lem}
\begin{proof}
    Clearly, the result is true when $W_{ij}=0$. Let $W_{ij} > 0$. Since the entries of $W$ are nonnegative it is easy to see that for every nonnegative integer $s$, $(W^s)_{ij} \geq 0.$ Moreover, for every positive integer $s$,
    $$
D^s =
\begin{cases}
\left(D_+ \right)^s-\left(D_- \right)^s, & \text{if } s \text{ is odd },\\
\left(D_+ \right)^s+\left(D_- \right)^s, & \text{if } s \text{ is even }.\\
\end{cases}
$$
\textbf{Case 1:} Let $p$ be an even integer. Then 
\begin{align}
   D^{p-1}&=\left(D_+ \right)^{p-1}-\left(D_- \right)^{p-1} \notag \\
   \implies \sigma W^{p-1} &= \left(D_+ \right)^{p-1}-\left(D_- \right)^{p-1}. \label{eq:p_even_sigmaW}
\end{align}
Using Lemma~\ref{lem:F_0for_posW} and \eqref{eq:p_even_sigmaW}, 
\begin{align*}
    (\left(D_+ \right)^{p-1})_{ij}&= \frac{\kappa}{\kappa+1} \sigma (W^{p-1})_{ij} \\
    \implies (\left(D_+ \right)^{p-1})_{ij}D_{ij}&= \frac{\kappa}{\kappa+1} \sigma^2 (W^{p-1})_{ij} W_{ij} \geq 0.
\end{align*}
\textbf{Case 2:} Let $p$ be an odd integer. Since 
\begin{equation*}
    0< \gamma -\delta= \operatorname{tr}(D^p)= \sigma \operatorname{tr}(W^p)
\end{equation*}
and $W^p$ is entrywise nonnegative, we have $\sigma=1$ and $D=W$. As $p-1$ is even we have
\begin{equation}\label{eq:p_odd}
    W^{p-1}=(D_+)^{p-1} + (D_-)^{p-1}
\end{equation}
Using Lemma~\ref{lem:F_0for_posW} and \eqref{eq:p_odd}, we get
\begin{align*}
    (\left(D_+ \right)^{p-1})_{ij} &=\frac{\kappa}{\kappa-1}(W^{p-1})_{ij} \geq 0 \\
    (\left(D_+ \right)^{p-1})_{ij} D_{ij}&=\frac{\kappa}{\kappa-1}(W^{p-1})_{ij} W_{ij}\geq 0.
\end{align*}
\textbf{Case 3:} Let $p \geq n-1$. \\
We have $p-1 \geq n-2$. Thus by Lemma~\ref{lem:DNNpower} $(W_+)^{p-1}$ is entrywise nonnegative. \\
If $\sigma=1$, then 
\begin{equation*}
    (\left(D_+ \right)^{p-1})_{ij} D_{ij}= (\left(W_+ \right)^{p-1})_{ij}  W_{ij} \geq 0.
\end{equation*}
If $\sigma=-1$ then $D_-=W_+$. By Lemma~\ref{lem:F_0for_posW},
\begin{equation*}
    (\left(D_+ \right)^{p-1})_{ij} =-\kappa (\left(W_+ \right)^{p-1})_{ij}
\end{equation*}
Therefore,
\begin{equation*}
    (\left(D_+ \right)^{p-1})_{ij}D_{ij} =\kappa (\left(W_+ \right)^{p-1})_{ij} W_{ij} \geq 0.
\end{equation*}

\end{proof}
 Now choose an arbitrary vector $k$-coloring of $G$. The coloring corresponds to a Gram matrix $Q$ with $Q \succeq 0, Q_{ii}=1  \quad (1 \leq i \leq n),~ Q_{ij} \leq \frac{-1}{k-1}$ for every edge $ij$. Fix an orthonormal eigenbasis $v^1,\ldots,v^n$ of $D$, with $Dv^i=d_i v^i$ for $1\le i\le n$. Define $X=(X_{ij})$ by 
$$X_{ij}= (v^i \circ v^j)^{\top}Q (v^i \circ v^j).$$
It is easy to see that
$$X_{ij} =\sum \limits_{l,t}Q_{lt}v^i_lv^j_lv^i_tv^j_t.$$

This is the construction used for \(P\) in the proof of Theorem~\ref{thm:chi_sv},
now applied to the orthonormal eigenbasis \(v^1,\ldots,v^n\) of \(D\). The same argument shows that \(X\) is symmetric and doubly stochastic. For further analysis, we define a function $\Phi(Z)=Q\circ Z.$ Note that $\Phi$ is a linear function.
\begin{lem}\label{lem:tr_formula}
    Let $f,g$ be real valued functions defined on the eigenvalues of $D$. Set $f_i=f(d_i), g_i=g(d_i)$. Then 
    $$\operatorname{tr}\left(f(D) \Phi\left(g(D)\right) \right)=\sum \limits_{i,j}X_{ij}f_ig_j.$$
\end{lem}
\begin{proof}
    \begin{align*}
        \operatorname{tr}\left(f(D) \Phi\left(g(D)\right) \right)&= \operatorname{tr}\left[ \left(\sum \limits_{i}f_iv^i(v^i)^\top\right)\Phi\left(\sum \limits_{j}g_jv^j(v^j)^\top \right) \right] \\
        &= \operatorname{tr}\left[ \left(\sum \limits_{i}f_iv^i(v^i)^\top\right)\left(\sum \limits_{j}g_j\Phi\left(v^j(v^j)^\top\right) \right) \right] \\
        &= \operatorname{tr} \left[ \sum \limits_{i} \sum \limits_{j} f_i g_j v^i (v^i)^\top \Phi(v^j(v^j)^\top) \right] \\
        &=\sum \limits_{i} \sum \limits_{j} f_i g_j \operatorname{tr}\left[ v^i (v^i)^\top \Phi(v^j(v^j)^\top)\right]\\
        &=\sum \limits_{i} \sum \limits_{j} f_i g_j \operatorname{tr}\left[ (v^i)^\top \Phi(v^j(v^j)^\top)v^i\right]\\
        &= \sum \limits_{i,j} f_i g_j \left( \sum_{l,t} Q_{lt}v^j_lv^j_tv^i_lv^i_t \right) \\
        &= \sum \limits_{i,j} f_i g_j \left[ (v^i \circ v^j)^\top Q (v^i \circ v^j) \right] \\
        &= \sum \limits_{i,j} X_{ij}f_i g_j.
    \end{align*}
\end{proof}
Recall that $d_i^+=\max\{d_i,0\}, d_i^-=\max\{-d_i, 0\}$. Then it can be shown that $z_1:= d_1^++\kappa^{\frac {1}{p-1}}d_1^-, \dots, z_n:= d_n^++\kappa^{\frac {1}{p-1}}d_n^-$ are the eigenvalues of $$H:= \left(D_+ \right)+\kappa^{\frac {1}{p-1}}\left(D_- \right)$$ with corresponding eigenvectors $v^1, \dots, v^n$. Since $d_i^+ d_i^-=0$, for every real $s > 0$ we have
$$z_i^s=(d_i^++\kappa^{\frac {1}{p-1}}d_i^-)^s=\left( d_i^+ \right)^s+(\kappa^{\frac {1}{p-1}})^s\left ( d_i^- \right)^s.$$
Therefore,
\begin{align*}
    H^s&= \sum \limits_{i} z_i^s v^{i} (v^{i})^{\top} \\
    &= \sum \limits_{i} \left(\left( d_i^+ \right)^s+(\kappa^{\frac {1}{p-1}})^s \left(d_i^- \right)^s\right)v^i(v^i)^\top \\
    &=\left(D_+ \right)^s+(\kappa^{\frac {1}{p-1}})^s\left(D_- \right)^s.
\end{align*}
In particular, 
\begin{equation}\label{eq:F_alt}
H^{p-1}=\left(D_+ \right)^{p-1}+\kappa \left(D_- \right)^{p-1}=F.
\end{equation}
\begin{lem}\label{lem:gam_lower}
    $$\operatorname{tr}\left(F\Phi(D_+)\right) \leq \gamma.$$
\end{lem}
\begin{proof}
    Consider
    \begin{align*}
        \operatorname{tr}\left(F\Phi(D)\right)&= \sum \limits_{ij} F_{ij} (\Phi(D))_{ji} \\
        &=\sum \limits_{ij} F_{ij} Q_{ij}D_{ij}.
    \end{align*}
    Using Lemma~\ref{lem:F_0for_posW}, whenever $D_{ij} \neq 0,$ we must have $F_{ij}=0$. Therefore we get 
    \begin{align*}
        \operatorname{tr}\left(F\Phi(D)\right)&=0\\
        \text{i.e. } \quad \operatorname{tr}\left(F\Phi(D_+)\right)&= \operatorname{tr}\left(F\Phi(D_-)\right).
    \end{align*}
    Therefore using \eqref{eq:F_alt},
    \begin{align*}
        \operatorname{tr}\left(H^{p-1}\Phi(H) \right)&= \operatorname{tr}\left(F\Phi(H)\right) \\
        &= \operatorname{tr}\left(F\Phi(D_+ )\right)+\kappa^{\frac {1}{p-1}}\operatorname{tr}\left(F\Phi(D_-)\right) \\
        &=\left(1+\kappa^{\frac {1}{p-1}}\right)\operatorname{tr}\left(F\Phi(D_+ )\right).
    \end{align*}
    Therefore,
    \begin{equation}\label{eq:tr_FPhiD+For} 
    \operatorname{tr}\left(F\Phi(D_+)\right)= \dfrac{\operatorname{tr}\left(H^{p-1}\Phi(H) \right)}{1+\kappa^{\frac{1}{p-1}}}.
    \end{equation}
    Now choose
    $$
g(t) =
\begin{cases}
t, & \text{if } t \geq 0,\\
\kappa^{\frac{1}{p-1}}(-t), & \text{if } t<0\\
\end{cases}
$$
and 
    $$
f(t) =
\begin{cases}
t^{p-1}, & \text{if } t \geq 0,\\
\kappa(-t)^{p-1}, & \text{if } t<0.\\
\end{cases}
$$
Clearly, $f(t)=g(t)^{p-1}$ for all $ t \in \mathbb{R}.$
Using Lemma~\ref{lem:tr_formula} we get,
\begin{equation}\label{eq:trH_p-1phi(H)}
    \operatorname{tr}(H^{p-1}\Phi(H))= \sum \limits_{i,j}X_{ij}z_i^{p-1}z_j.
\end{equation}
Now consider,
\begin{equation*}
    \sum \limits_{i,j} X_{ij}z_i^{p}= \sum \limits_{i}z_i^{p} \sum \limits_{j}X_{ij}=\sum \limits_{i}z_i^{p}.
\end{equation*}
Similarly, 
\begin{equation*}
    \sum \limits_{i,j} X_{ij}z_j^{p}= \sum \limits_{j}z_j^{p}.
\end{equation*}
Therefore,
\begin{equation}\label{eq:tr_Hp}
    \sum \limits_{i}z_i^{p}=\frac 12 \sum \limits_{i,j} X_{ij}\left(z_i^{p}+z_j^{p}\right).
\end{equation}
Also by symmetry of $X$,
\begin{equation*}
    \sum \limits_{i,j}X_{ij}z_i^{p-1}z_j = \sum \limits_{i,j}X_{ij}z_j^{p-1}z_i.
\end{equation*}
Therefore,
\begin{equation}\label{eq:trH_p-1Phi(H)For}
    \sum \limits_{i,j}X_{ij}z_i^{p-1}z_j= \frac 12 \sum \limits_{i,j}X_{ij}\left(z_i^{p-1}z_j+z_j^{p-1}z_i \right).
\end{equation}
Using \eqref{eq:trH_p-1phi(H)}, \eqref{eq:tr_Hp} and \eqref{eq:trH_p-1Phi(H)For} we get,
\begin{align*}
    \operatorname{tr}(H^{p})-\operatorname{tr}(H^{p-1}\Phi(H))&= \frac 12 \sum \limits_{i,j}X_{ij}\left(z_i^p+z_j^p-z_i^{p-1}z_j-z_j^{p-1}z_i \right) \\
    &=\frac 12 \sum \limits_{i,j}X_{ij}\left(z_i^{p-1}-z_j^{p-1}\right)\left(z_i-z_j\right).
\end{align*}
Since $p-1>0$ we get,
$$\operatorname{tr}(H^p)\geq \operatorname{tr}(H^{p-1}\Phi(H)).$$
Then \eqref{eq:tr_FPhiD+For} becomes,
\begin{equation*}
    \operatorname{tr}(F\Phi(D_+))\leq \dfrac{\operatorname{tr}(H^p)}{1+\kappa^{\frac{1}{p-1}}}.
\end{equation*}
Finally, note that
\begin{align*}
    \operatorname{tr}(H^p)&=\operatorname{tr}(\left(D_+ \right)^p+\kappa^{\frac{p}{p-1}}\left(D_- \right)^p) \\
    &=\operatorname{tr}(\left(D_+ \right)^p)+\kappa^{\frac{p}{p-1}}\operatorname{tr}(\left(D_- \right)^p)\\
    &=\gamma+\kappa^{\frac{p}{p-1}}\delta\\&=\gamma+\kappa^{\frac{1}{p-1}}\kappa \delta \\
    &=(1+\kappa^{\frac{1}{p-1}})\gamma.
\end{align*}
Therefore we get,
$$\operatorname{tr}(F\Phi(D_+))\leq \gamma.$$
\end{proof}
We now use following notations for further calculations.
$$L= \operatorname{tr}\left(\left(D_+ \right)^{p-1}\Phi(D_-)\right) \quad \quad M=\operatorname{tr}\left(\left(D_- \right)^{p-1}\Phi(D_+ )\right).$$ In the next lemma, we derive an upper bound for $M$.
\begin{lem}\label{lem:M_upper}
For $M$ as defined above, we have the following:
    $$M \leq \delta.$$
\end{lem}
\begin{proof}
Using $g(t)=t^+$ and $f(t)=(t^+)^{p-1}$ with Lemma~\ref{lem:tr_formula} we get
  \begin{equation}\label{eq:tr_D+_nonneg}
  \operatorname{tr}\left(\left(D_+ \right)^{p-1}\Phi(D_+ )\right)= \sum \limits_{ij} X_{ij}(d_i^+)^{p-1}d_j^+ \geq 0.
  \end{equation}
  Therefore,
  \begin{align}\label{eq:kapM_upper}
      \operatorname{tr}\left(F\Phi(D_+ )\right)&=\operatorname{tr}\left( \left(\left(D_+ \right)^{p-1}+\kappa \left(D_- \right)^{p-1}\right) \Phi(D_+ ) \right)  \notag \\
      &=\operatorname{tr}\left(\left(D_+ \right)^{p-1}\Phi(D_+ )\right)+\kappa \operatorname{tr}\left(\left(D_- \right)^{p-1}\Phi(D_+)\right)  \notag \\
      &\geq \kappa M. 
  \end{align}
 Using Lemma~\ref{lem:gam_lower} and \eqref{eq:kapM_upper} we get
  $\kappa M \leq \gamma$. Therefore,
  $$M \leq \dfrac{\gamma}{\kappa} =\delta.$$
\end{proof}
In the next lemma, we derive lower and upper bounds for $L$.
\begin{lem}\label{lem:L_upperlower}
For $L$ and $M$ as defined before we have
    $$\dfrac{1}{k-1}\gamma \leq L \leq M^{\frac {1}{p-1}}\gamma^{\frac{p-2}{p-1}}.$$
\end{lem}
\begin{proof}
    If $D_{ij} \neq 0$, then $\{i, j\}$ is an edge and hence $Q_{ij} \leq -\frac{1}{k-1}$. Therefore using Lemma~\ref{lem:D_+D_geq0} we obtain
    $$(\left(D_+ \right)^{p-1})_{ij}D_{ij}Q_{ij} \leq -\frac{1}{k-1} (\left(D_+ \right)^{p-1})_{ij}D_{ij}.$$
    If $D_{ij}=0$, the above inequality is again true. Therefore,
    \begin{align*}
        \operatorname{tr}(\left(D_+ \right)^{p-1}\Phi(D)) &= \sum \limits_{i,j} (\left(D_+ \right)^{p-1})_{ij}D_{ij}Q_{ij} \\
        &\leq -\frac{1}{k-1} \sum \limits_{i,j} (\left(D_+ \right)^{p-1})_{ij}D_{ij}= -\frac{1}{k-1}\operatorname{tr}(\left(D_+ \right)^{p-1}D) \\
        &= -\frac{1}{k-1}\operatorname{tr}(\left(D_+ \right)^{p-1}(D_+ - D_-)) \\
        &= -\frac{1}{k-1}\operatorname{tr}(\left(D_+ \right)^{p}) =-\frac{1}{k-1}\gamma.
    \end{align*}
    Using \eqref{eq:tr_D+_nonneg},
    \begin{align*}
        \operatorname{tr}(\left(D_+ \right)^{p-1}\Phi(D))&=\operatorname{tr}(\left(D_+ \right)^{p-1}\Phi(D_+))-\operatorname{tr}(\left(D_+ \right)^{p-1}\Phi(D_- )) \\
        &\geq -L.
    \end{align*}
    Combining the above two inequalities we get, 
    $$L \geq \dfrac{1}{k-1}\gamma.$$
    This proves the lower bound.
    
    Next we derive the upper bound. Let $g(t)=t^-$ and $f(t)=(t^+)^{p-1}$. Then, by Lemma~\ref{lem:tr_formula}, we get
    \begin{equation*}
        L=\sum \limits_{i,j} X_{ij}\left(d_i^+\right)^{p-1}d_j^-.
    \end{equation*}
    Similarly for the choice $g(t)=t^+$ and $ f(t)=\left( t^- \right)^{p-1}$, we get
    \begin{equation*}
        M=\sum \limits_{i,j} X_{ij}\left(d_i^-\right)^{p-1}d_j^+.
    \end{equation*}
    Using the symmetry of $X$ we get
    \begin{equation}\label{eq:M_for}
        M=\sum \limits_{i,j} X_{ji}\left(d_i^-\right)^{p-1}d_j^+= \sum \limits_{i,j} X_{ij}d_i^+\left(d_j^-\right)^{p-1}.
    \end{equation}
    Note that,
    \begin{equation}\label{eq:gam_alter}
        \sum \limits_{i,j} X_{ij}\left(d_i^+\right)^{p} = \sum \limits_{i} \left(d_i^+\right)^{p}\sum \limits_{j} X_{ij}=\gamma.
    \end{equation}
    Suppose $p>2$. Note that,
    $$X_{ij}\left(d_i^+\right)^{p-1}d_j^- = \left( X_{ij}d_i^+ \left(d_j^-\right)^{p-1} \right )^{\frac{1}{p-1}} \left( X_{ij}\left(d_i^+\right)^{p}\right)^{\frac{p-2}{p-1}}.$$
    Using the H\"older's inequality with conjugate exponents $p-1$ and $(p-1)/(p-2)$, we get,
    \begin{align*}
        L &\leq \left(\sum \limits_{i,j} X_{ij}d_i^+ \left( d_j^-\right)^{p-1} \right )^{\frac{1}{p-1}} \left(\sum \limits_{i,j}  X_{ij}\left (d_i^+\right)^{p}\right)^{\frac{p-2}{p-1}} \\
        &= M^{\frac{1}{p-1}}\gamma^{\frac{p-2}{p-1}}.
    \end{align*}
    where the last equality comes using \eqref{eq:M_for} and \eqref{eq:gam_alter}.\\
    If $p=2$, then $L=M$ and hence the upper bound clearly holds.
    
\end{proof}
\noindent We now complete the proof of Theorem~\ref{thm:partial}.\\
\begin{proof}
Using Lemma~\ref{lem:M_upper} and Lemma~\ref{lem:L_upperlower} we get,
$$\dfrac{1}{k-1}\leq \dfrac{L}{\gamma}\leq \left( \dfrac{M}{\gamma}\right)^{\frac{1}{p-1}}\leq \left( \dfrac{\delta}{\gamma}\right)^{\frac{1}{p-1}}=\kappa^{-\frac{1}{p-1}}.$$
Therefore, we get
\begin{equation}\label{eq:kap_upper}
    \kappa^{\frac{1}{p-1}}\leq k-1.
\end{equation}

Since $A/m\in\mathcal W$, homogeneity gives $\kappa_A(p)=\kappa_{A/m}(p)\leq \kappa$.
Thus using \eqref{eq:kap_upper} we get
$$\left(\kappa_A(p) \right)^{\frac{1}{p-1}}\leq k-1.$$
Hence
$$k \geq 1+\left(\kappa_A(p) \right)^{\frac{1}{p-1}}.$$
The RHS is a lower bound for $k$ if $G$ has a proper vector $k-$coloring. But $\chi_v(G)$ is the infimum of all such feasible values of $k$. Therefore, 
$$\chi_v(G) \geq 1+ \left(\max \left \{ \dfrac{\operatorname{tr}\left((A_+)^p\right)}{\operatorname{tr}\left((A_-)^p\right)}, \dfrac{\operatorname{tr}\left((A_-)^p\right)}{\operatorname{tr}\left((A_+)^p\right)}\right \}\right)^{\frac{1}{p-1}}.$$
\end{proof}
\begin{rem}
    From the proof of Theorem~\ref{thm:partial}, it is easy to see that Theorem~\ref{thm:partial} holds with $A$ replaced by any nonzero entrywise nonnegative weighted Hermitian adjacency matrix $B=(B_{uv})$ of $G$ since $B$ is real symmetric and for $s=\sum \limits_{uv \in E(G)} B_{uv}$ we have $B/s \in \mathcal{W}$.
\end{rem}
\section{Counterexamples}\label{sec:counter}
We first consider the extension of the bound in Conjecture~\ref{con:chi_v_conj} to $p<1.5$. Elphick, Tang and Zhang \cite{MR4971662} showed that $C_5$ is a counterexample for the $p$-energy bound being a lower bound for $\chi_v$ for $0 \le p \le 1.154$ with exponent $1/|p-1|$ for \(p<1\) and the limiting
interpretation at $p=1$.

All connected Strongly Regular Graphs (SRG)s have $\chi_v(G) = 1 + \lambda_1(G)/|\lambda_n(G)|$ (see Lemma~\ref{lem:chi_vSRG}), which makes it easier to search for  counterexamples. We extend the list of counterexamples for values of $p$ in the range $(1,3/2)$. 

A finite generalised quadrangle is a point--line incidence
structure whose incidence graph is connected and bipartite,
with diameter \(4\) and girth \(8\). It has order \((s,t)\)
if each line contains \(s+1\) points and each point lies on
\(t+1\) lines. Its collinearity graph has the points as vertices,
with two distinct points adjacent precisely when they lie on a
common line \cite[Section~9.6]{BrouwerHaemersBook}.

We first recall a few known facts about the generalised quadrangles. For every prime power \(q\), there exists a generalised quadrangle of order $(q,q^2)$. Fix such a quadrangle and let \(G_q\) denote its collinearity graph. Elphick and Linz~\cite[proof of Theorem~3.1]{MR4755798} used this same family of collinearity graphs to show that the ratio of the negative $2$-energy to the positive $2$-energy can grow as \(\Theta\bigl(|V(G_q)|^{1/4}\bigr)\). It is known that $G_q$ is strongly regular with parameters $\left( (q+1)(q^3+1), q(q^2+1),q-1,q^2+1 \right).$ The spectrum of $G_q$ is $\{ q(q^2+1)^{(1)}, (q-1)^{(q^4+q^2)}, -(q^2+1)^{(q^3-q^2+q)}\}$ (see \cite[Section~9.6]{BrouwerHaemersBook}). Also the following is a known result.
\begin{lem}[{\cite[Lemma~5.2]{GodsilRobersonRobertSimone2016}}]\label{lem:chi_vSRG}
    If $G$ is a connected SRG,
    $$\chi_v(G)=\chi_{sv}(G)=1-\dfrac {\lambda_1(G)}{\lambda_n(G)}.$$
\end{lem}
Therefore we get
\begin{equation*}
    \chi_v(G_q)=\chi_{sv}(G_q)= 1+\dfrac{q(q^2+1)}{(q^2+1)}=q+1.
\end{equation*}
We are now ready to provide the proposed counterexamples as follows.
\begin{thm}\label{thm:smallrange_counter}
    Fix $1<p<3/2$. Let $q \geq 2$  be a prime power satisfying $q \geq (p-1)^\frac{-1}{3-2p}$. Then 
    $$1+\max\left \{ \dfrac{\mathcal{E}_p^+(G_q)}{\mathcal{E}_p^-(G_q)}, \dfrac{\mathcal{E}_p^-(G_q)}{\mathcal{E}_p^+(G_q)}\right \}^{\dfrac {1}{p-1}}> \chi_v(G_q)=\chi_{sv}(G_q)=q+1.$$
\end{thm}
\begin{proof}
    Since 
    $$0 < \dfrac{q(q-1)}{q^2+1}< \dfrac{q-1}{q}=1-\dfrac{1}{q}$$
    and $0 < p-1 < 1$, we obtain 
    \begin{equation}\label{eq:q_inter_lower}
        \left( \dfrac{q(q-1)}{q^2+1} \right)^{p-1} \leq \left( 1-\dfrac{1}{q} \right)^{p-1}.
    \end{equation}
    Using weighted AM-GM inequality we obtain,
    \begin{align}\label{eq:q_inter_upper}
        \left( 1-\dfrac{1}{q} \right)^{p-1} &\leq (p-1) \left( 1-\dfrac{1}{q} \right)+(2-p)\cdot1  \notag \\
        &=1-\dfrac{p-1}{q}.
    \end{align}
    Combining \eqref{eq:q_inter_lower} and \eqref{eq:q_inter_upper} we obtain 
    \begin{equation}\label{eq:q_used_upper}
        \left( \dfrac{q(q-1)}{q^2+1} \right)^{p-1} \leq 1-\dfrac{p-1}{q}.
    \end{equation}
    Moreover, the condition on $q$ is equivalent to 
    \begin{equation}\label{eq:q_alt}
        q^{2p-1} \leq (p-1)q^2.
    \end{equation}
    Therefore using \eqref{eq:q_alt} and \eqref{eq:q_used_upper} we get
    \begin{align*}
        \dfrac{q^{p-1} \mathcal{E}_p^+(G_q)}{(q^2+1)^p} &= q^{2p-1} + \dfrac{q^{p+1}(q-1)^p}{(q^2+1)^{p-1}}\\
        &= q^{2p-1} +q^2(q-1) \left( \dfrac{q(q-1)}{q^2+1} \right)^{p-1} \\
        &\leq (p-1)q^2 +q^2 (q-1)\left( 1-\dfrac{p-1}{q}\right)\\
        &= q^3 -q^2 +(p-1)q \\
        &< q^3 -q^2 +q \\
        &= \dfrac{\mathcal{E}_p^-(G_q)}{(q^2+1)^p}.
    \end{align*}
    where the strict inequality holds because $p-1 < 1$.  Therefore we obtain
    $$\mathcal{E}_p^-(G_q) > q^{p-1} \mathcal{E}_p^+(G_q).$$
    Since $q^{p-1}>1$ we get
    $$ 1+\max\left \{ \dfrac{\mathcal{E}_p^+(G_q)}{\mathcal{E}_p^-(G_q)}, \dfrac{\mathcal{E}_p^-(G_q)}{\mathcal{E}_p^+(G_q)}\right \}^{\dfrac{1}{p-1}}=1+ \left( \dfrac{\mathcal{E}_p^-(G_q)}{\mathcal{E}_p^+(G_q)}\right)^{\dfrac{1}{p-1}}> q+1 =\chi_{sv}(G_q)=\chi_v(G_q).$$
\end{proof}

Since powers of $2$ are unbounded, an admissible prime power $q$ exists for every fixed $1<p<3/2$. Moreover, $C_5$ is strongly regular, so Lemma~\ref{lem:chi_vSRG} gives $\chi_v(C_5)=\chi_{sv}(C_5)=\sqrt5$. Combining Theorem~\ref{thm:smallrange_counter} with the $C_5$ example of Elphick, Tang and Zhang~\cite{MR4971662}, we obtain, for every $0\le p<3/2$, a graph for which the extended bound fails even for $\chi_{sv}$.

\begin{rem}
    With the assistance of ChatGPT (OpenAI), we identified counterexamples for each fixed non integer $p>2$ to the conjecture of Elphick, Tang and Zhang~\cite[Conjecture~6.3]{MR4971662}, which will be presented in a follow-up manuscript, as the construction is involved. 
\end{rem}

\section{Conclusion}\label{sec:conc}

Conjecture 6.3 in Elphick, Tang and Zhang \cite{MR4971662}  asks whether for any graph on $n$ vertices:

\[
\chi_v(G) \ge 1 +\max\left \{ \frac{\mathcal{E}_p^+(G)}{\mathcal{E}_p^-(G)}, \frac{\mathcal{E}_p^-(G)}{\mathcal{E}_p^+(G)} \right \}^{\frac{1}{p - 1}} \mbox{  for  } p  \ge 2?
\]
We have proved the analogous bound for the strict vector chromatic number. We also have proved the conjecture for all integer values of $p \geq 2$ and for all reals $p \geq n-1$. We  present SRGs (Generalised Quadrangles) for which this bound is false for $1< p < 1.5$, which improves on the limit of $ p = 1.154 $ in \cite{MR4971662}. The conjectured bound for the vector chromatic number is false in general. Indeed, for all non-integral $p \ge 2$ there exists a counter-example which we shall present in a follow-up manuscript. 

\section*{Acknowledgements}
The authors acknowledge Dr. Shivaramakrishna Pragada for helpful discussions. Abhay Jayarajan acknowledges support from the DST INSPIRE Fellowship (IF220692). M.~Rajesh Kannan acknowledges financial support from ANRF-CRG, India (File No.\ CRG/2023/002747). 	

\section*{Declaration of AI use}

The authors acknowledge use of ChatGPT (GPT-6 Astra, OpenAI; accessed October and November 2026) to find the family of GQ graphs and the counterexamples for the conjecture in general, to check for typographical errors and to refine some of the proofs. All mathematical statements have been independently verified by the authors, who are solely responsible for the correctness of the results presented in this paper. 

\bibliographystyle{plain}
\bibliography{Ref}

\vspace{2cm}

\affl{Clive Elphick}{clive.elphick@gmail.com}{School of Mathematics, University of Birmingham, Birmingham, UK}

\affl{Abhay Jayarajan}{ma23resch02001@iith.ac.in, abhayjayarajan@gmail.com}{ Department of Mathematics, Indian Institute of Technology Hyderabad, Kandi, Sangareddy 502284, India.}

\affl{M. Rajesh Kannan}{rajeshkannan@math.iith.ac.in, rajeshkannan1.m@gmail.com}{Department of Mathematics, Indian Institute of Technology Hyderabad, Sangareddy 502285, India}

\end{document}